\documentclass[12pt]{article}
\usepackage{amsmath}
\usepackage{amsthm}
\usepackage{amssymb}
\usepackage{mathrsfs}
\usepackage{authblk}
\usepackage[usenames]{color}
\usepackage[latin5]{inputenc}
\usepackage{cite}
\usepackage{pdfsync}
\usepackage{enumitem}
\setlist{parsep=0pt,itemindent=0pt}

\usepackage{graphicx}

\theoremstyle{definition}

\newtheorem{thm}{Theorem}

\numberwithin{equation}{section}
\numberwithin{thm}{section}
\numberwithin{lemma}{section}
\numberwithin{prop}{section}
\numberwithin{cor}{section}
\numberwithin{rmk}{section}
\numberwithin{defn}{section}
\numberwithin{exa}{section}

\def\frac#1#2{{#1\over #2}}

\begin{document}
\pagenumbering{arabic}
\clearpage
\thispagestyle{empty}

\title{Integrability of the Basener-Ross model with time-dependent
coefficients}

\author[1]{F.~G\"ung\"or\thanks{gungorf@itu.edu.tr}}
\author[2]{P.J. Torres\thanks{ptorres@ugr.es}}

\affil[1]{Department of Mathematics, Faculty of Science and
Letters, Istanbul Technical University, 34469 Istanbul, Turkey}
\affil[2]{Departamento de Matem\'atica Aplicada and Excellent Research Unit ``Modeling Nature'' (MNat), Universidad de
Granada, 18071 Granada, Spain}

\date{\today}

\maketitle

\begin{abstract} The Basener-Ross system is a known model in Population Dynamics for the interaction of consumers and resources in an isolated habitat.
For an extended version with time-dependent coefficients as a model
of possible variations of the environtmental conditions, some relations among the coefficients are provided leading to the integrability of the system.
\end{abstract}

\vspace{1,5mm}




\section{Introduction}

In Ecology, an isolated island is an excellent laboratory to explore the evolution of ecosystems due to the total absence of external distortion factors like migration. In the paper \cite{BasenerRoss2004}, Basener and Ross proposed the following model for the evolution of human population in Easter island
\begin{equation}\label{model_constant}
\begin{array}{l}
\displaystyle{\dot{x} = c x\left(1-\frac{x}{\kappa}\right)-h y},\\
\\
\displaystyle{\dot{y} = a y\left(1-\frac{y}{x}\right)},
\end{array}
\end{equation}
where $a ,c,h,\kappa$  are positive constants, $x(t), y(t)$ are
respectively the amount of resources and the human population in
the island at time $t$. The interaction is of predator-prey type,
in such a way that the amount of resources $x$ serves as a dynamic
carrying capacity for the predator $y$. Comparatively with other
predator-prey systems available in the literature, the
Basener-Ross model presents a remarkable variety of dynamical
behaviours, including extinction at a finite time. By this reason,
it has been recognized as a plausible model for evolution of
population in ancient civilizations and several generalizations
and variants have been formulated and studied, see for instance
\cite{BB,BF,BM,K} and the references therein.

A reasonable extension of the Basener-Ross model is the
consideration of time-changing ecological parameters. Variation in
time of any of the parameters involved in system
\eqref{model_constant} can be easily justified from a biological
point of view. For example, climate fluctuations may influence
critically on the dynamics of interacting species by changing the
natural growth rates (see e.g. \cite{Am,ZT,FF,HB,SR} and the
references therein). Also, empirical evidence shows how
technological advances may increase the carrying capacity $\kappa$
of a resource \cite{C,MA}. Finally, time-variation in $h$ models
the adoption of a suitable harvesting protocol.

According to the previous discussion, we are  going to study the non-autonomous system
\begin{equation}\label{model_eq1}
\begin{array}{l}
\displaystyle{\dot{x} = c(t) x\left(1-\frac{x}{\kappa(t)}\right)-h(t) y},\\
\\
\displaystyle{\dot{y} = a(t) y\left(1-\frac{y}{x}\right)},
\end{array}
\end{equation}
where now $a(t),c(t),\kappa(t),h(t)$ are continuous functions with positive values. Our main motivation comes from the recent paper \cite{NucciSanchini2015}, where the integrability of system \eqref{model_constant} is studied. Our objective is to identify novel families of integrable cases for the system with variable coefficients \eqref{model_eq1} that are new even in the case of constant coefficients. The main tool will be the reduction to a second order differential equation and the application of recent integrability criteria developed in \cite{CGT}.

\section{Reduction to a second order ODE and integrable cases}

As a first step, by means of the change of time variable
\begin{equation}\label{change}\tau=\int_0^t a(s)ds,
\end{equation}
we can consider without any loss of generality that $a(t)\equiv 1$. With this change, the system reads
\begin{equation}\label{model_eq2}
\begin{array}{l}
\displaystyle{\dot{x} = c(t) x\left(1-\frac{x}{\kappa(t)}\right)-h(t) y},\\
\\
\displaystyle{\dot{y} = y\left(1-\frac{y}{x}\right)}.
\end{array}
\end{equation}

From the second equation, one finds
\begin{equation}\label{x}
 x=\frac{y^2}{y-\dot{y}}.
\end{equation}
Then, elimination of $x$ in \eqref{model_eq2} leads to the second order ODE
\begin{equation}\label{2nd-model-1}
  \ddot{y}+r(t)\dot{y}+\alpha(t)y=\beta(t)\frac{\dot{y}^2}{y}+\gamma(t)y^2,
\end{equation}
where the coefficients are defined by
\begin{equation}\label{coeff-model-1}
  r(t)=1+c(t)-2h(t),  \quad \alpha(t)=h(t)-c(t),  \quad \beta(t)=2-h(t), \quad \gamma(t)=-\frac{c(t)}{\kappa(t)}.
\end{equation}

Now, the objective is to show that \eqref{2nd-model-1} is  integrable under some particular relations of its coefficients. A first equation that is solvable by quadratures is obtained by restricting the arbitrary function $H$ in \cite[(3.13)]{CGT} to $H(I)=4AI^2-B$, leading to
\begin{equation}\label{solvable}
 \ddot{y}-A\nu\dot{y}+\left[2\dot{\nu}+(1-A)\nu^2\right]y=\frac{1}{4}\left(5+A\right)\frac{\dot{y}^2}{y}-By^2,
\end{equation}
where the constants $A,B$ are arbitrary and the function $\nu(t)$ can be chosen in the form $\nu=\dot{a}/a$ with $a(t)$ being any solution of the equation
\begin{equation}\label{3rd}
  \dddot{a}+4p\dot{a}+2\dot{p}a=0
\end{equation}
for a given $p(t)$. Eliminating $a$ gives a second order ODE for  $\nu$
\begin{equation}\label{nu-p-eq}
  \ddot{\nu}+3\nu \dot{\nu}+\nu^3+4p(t)\nu+2\dot{p}(t)=0
\end{equation}
This equation belongs to the particular form of the second member of the
Riccati chain (see \cite{GL} for details). The substitution
$\nu=\dot{a}/a$ plays the role of a Cole-Hopf transformation  that
achieves linearizability of the latter equation.

A second equation of interest is the following dissipative form of
Kummer-Schwarz equation
\begin{equation}\label{diss-KS}
  \ddot{w}+r(t)\dot{w}+\frac{4p(t)}{1-n}w=\sigma\frac{\dot{w}^2}{w}+\frac{4q}{1-n}\exp\left[-2\int_0^t r(\tau)d\tau\right]w^n, \quad n\ne 1,
  \quad \sigma=\frac{n+3}{4},
\end{equation}
that has been introduced in \cite[(3.37)]{CGT}. The linear transformation
\begin{equation}\label{trans}
w=\phi(t)z(t), \quad
\phi(t)=\exp\left[\frac{1}{2(\sigma-1)}\int_0^t r(\tau)d\tau\right],
\end{equation}
transforms \eqref{diss-KS} to
\begin{equation}\label{stand}
  \ddot{z}+\frac{4}{1-n}I(t)z=\sigma\frac{\dot{z}^2}{z}+\frac{4q}{1-n}z^n,
\end{equation}
where $$I(t)=p-\frac{1}{4}(r^2+2\dot{r}).$$
It is known that Eq. \eqref{stand} has the general solution
\begin{equation}\label{super}
  z=(Au^2+2Buv+Cv^2)^{2/(1-n)},  \quad (AC-B^2)W^2=q,
\end{equation}
where $u, v$ are two linearly independent solutions of the base equation
\begin{equation}\label{base}
  \ddot{z}+I(t)z=\ddot{z}+[p-\frac{1}{4}(r^2+2\dot{r})]z=0.
\end{equation}
In consequence, the general solution of  \eqref{diss-KS} is given by
\begin{equation}\label{gen-sol-dissip}
 w=\exp\left[\frac{2}{n-1}\int_0^t r(\tau)d\tau\right](Au^2+2Buv+Cv^2)^{2/(1-n)},  \quad (AC-B^2)W^2=q.
\end{equation}

\section{A first family of solvable models}

We want to analyze relations among the coefficients leading to the
integrable cases. First of all, Eq. \eqref{2nd-model-1} is
linearizable by the transformation $y=1/z$ only when
$h(t)\equiv0$. For $h\ne 0$,  a family of integrable systems can
be found by adapting \eqref{2nd-model-1} to the ansantz
\eqref{solvable}. By identifying the coefficients of both
equations, we get
\begin{equation}\label{link1}
c(t)=B\kappa(t),\quad h(t)=\frac{3-A}{4}
\end{equation}
and
\begin{equation}\label{link2}
h(t)-c(t)=2\dot\nu+\left(1-A\right)\nu^2,\quad 1+c(t)-2h(t)=-A\nu.
\end{equation}
From this information, $h$ must be constant, and then $c(t)$ must verify  simultaneously two conditions
$$
c(t)=h-2\dot\nu+2\left(1-2h\right)\nu^2,\quad c(t)=2h-1-(3-4h)\nu.
$$
Hence, $\nu$ must be a solution of the equation
$$
h-2\dot\nu+2\left(1-2h\right)\nu^2=2h-1-(3-4h)\nu
$$
with $h$ an arbitrary constant, or equivalently
\begin{equation}\label{nu}
2\dot\nu=1-h+(3-4h)\nu+\left(2-4h\right)\nu^2.
\end{equation}
Moreover, this condition must be satisfied together with \eqref{nu-p-eq}. Substitution of \eqref{nu} into \eqref{nu-p-eq} gives the finite constraint (a cubic equation for $\nu$ or a quadratic equation for $h$)
\begin{equation}\label{cond}
       4(4h-3)(h-1)\nu^2(2\nu+3)
      +[6(4h-3)(h-1)+1+16p)]\nu+(4h-3)(h-1)+8\dot{p}=0.
\end{equation}
For constant choice of $p(t)$, $\nu$ has to be constant. For instance, taking $p=-\nu^2/4$ ($\nu$ constant) we find that the following two values of $\nu$ as solutions of \eqref{nu} also satisfy \eqref{cond}
\begin{equation}\label{nu-h}
  \nu=\frac{1-h}{2h-1},  \quad \nu=-\frac{1}{2}.
\end{equation}

Summing up, the following result holds.

\begin{thm}
For a fixed $h>0$, let $\nu(t)$ be a solution of the system \eqref{nu}-\eqref{cond} for some $p(t)$. Then, system \eqref{model_eq2} is solvable for the coefficients
$$
h(t)\equiv h,\quad c(t)=2h-1-(3-4h)\nu(t),\quad  \kappa(t)=\frac{c(t)}{B},
$$
where $B$ is an arbitrary constant.
\end{thm}

The case of constant coefficients comes from taking $\nu$ as a
constant in conditions \eqref{link2} so that from \eqref{nu-h} we have either $\nu=-1/2$ ($c=1/2$) or
\begin{equation}\label{solpar}
h=\frac{1+\nu}{1+2\nu},\quad c=1-\nu,
\end{equation}
that can be seen as a curve in parametric coordinates for which
the model is solvable. Eliminating the parameter we arrive at the
explicit curve
\begin{equation}\label{curve}
h(c)=\frac{2-c}{3-2c}.
\end{equation}
Of course, to be biologically relevant both parameters $h,c$ must
be positive. Figure 1 is a graphical representation of the curve.
It coincides with the curve identified in \cite{NucciSanchini2015}
(see (29) therein).

\begin{figure}\label{fig1}
\begin{center}\includegraphics{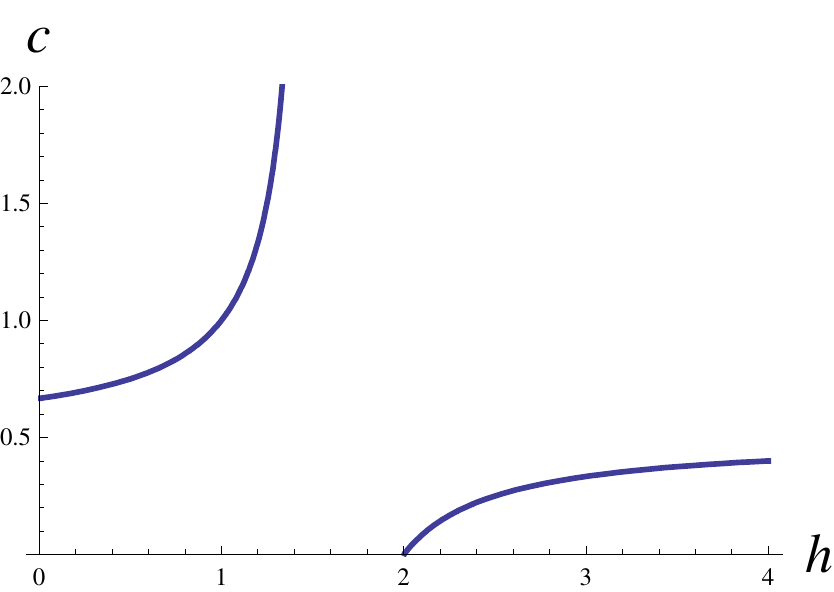}
\end{center}
\caption{The curve of admissible values in the plane $h-c$.}
\end{figure}

Under the condition \eqref{curve}, we know from \cite{CGT} that
the change of variables
\begin{equation}\label{change-bio}
r=e^{-\frac{\nu}{2} t} y,  \quad s=-\frac{1}{\nu}e^{-\nu t}
\end{equation}
transforms  \eqref{solvable} into the equivalent equation
\begin{equation}\label{eq-r}
r''(s)=\frac{5+A}{4}\frac{r'^2}{r}-B r^2.
\end{equation}
As shown in \cite{CGT}, this equation can be reduced to
quadratures. To this end, we let $R=dr/ds$ and exchange the roles
of $(r,s)$. Using the coordinates
$$\omega=-\frac{1}{4}r^{-3/2}R,  \quad \xi=\ln r$$ we find the equation
\begin{equation}\label{sep}
 \frac{d\omega}{d\xi}=\frac{4(A-1)\omega^2-B}{16\omega}.
\end{equation}
This equation is separable and of Bernoulli type, then it is easy
to find the explicit solution $\omega=\Phi(\ln r,c_1)$, where $c_1$ is an arbitrary constant. The
general solution is obtained by integrating the separable first
order ODE
\begin{equation}\label{genA}
\frac{dr}{ds}=-4r^{3/2}\Phi(\ln r,c_1).
\end{equation}
At this moment, it is necessary to analyze the cases $A>1$ and $A<1$ separately. For $A> 1$, some simple computations lead to the solution by quadratures of \eqref{genA} as
\begin{equation}\label{int}
-\frac{\sqrt{A-1}}{2}\int^{r(s)}_0 \frac{r^{-3/2}dr}{\sqrt{c_1 r^{(A-1)/2}+B}}=s+c_2.
\end{equation}
For general values $A>1$ this integral can be expressed in terms of hypergeometric functions as
\begin{equation}\label{HGF}
  \frac{\sqrt{2\delta-1}}{\sqrt{2c_1}\delta}r^{\delta/(1-2\delta)}
  {}_2F_1\left(\delta,\frac{1}{2};\delta+1;-\frac{B}{c_1}r^{1/(2\delta-1)}\right),
\end{equation}
where
$\delta=\frac{A+1}{2(A-1)}>1/2.$

There are many values of $\delta$ (or $A$) for which the
hypergeometric function is elementary. For example, in the
particular case $A=2$, this last integral is solvable
and $r(s)$ can be written as
\begin{equation}\label{int}
r(s)=\frac{B}{B^2(s+c_2)^2-c_1},
\end{equation}
where $c_1>0$ and $c_2$ is an arbitrary constant. Note that $A=2$ means $h=1/4$ and $\nu= -3/2$ by \eqref{nu-h}.  Then, by inverting the change \eqref{change-bio},
we find the general solution of \eqref{solvable} for $A=2$ as
$$y(t)=B\left(\frac{4}{9}B^2 e^{3 t}+(B^2c_2^2-c_1)e^{-3 t}+\frac{4}{3}B^2c_2\right)^{-1}.$$
As a consequence, the general solution of the Basener-Ross system \eqref{model_eq2} is readily obtained by using \eqref{x}.

The case $A<1$ can be analyzed as before. In this case, the analogous equation to \eqref{int} is now
\begin{equation}\label{int2}
-\frac{\sqrt{1-A}}{2}\int^{r(s)}_0 \frac{r^{-3/2}dr}{\sqrt{c_1 r^{(A-1)/2}-B}}=s+c_2.
\end{equation}
As before, in general the solution will be obtained in an implicit form by quadratures. The case $A=0$ admits an explicit form of the solution. Since the analysis is similar, we skip further details.

\section{A second family of solvable models}

In this section, we are going to use the dissipative KS equation
\eqref{diss-KS} for the identification of a new class of
coefficients enabling the integrability of Eq.
\eqref{2nd-model-1}. By comparing \eqref{2nd-model-1} both
equations, we find $n=2$ and the relations
\begin{equation}\label{relat-1}
 4q \exp\left[-2\int_0^t (1+c(\tau)-2h(\tau))d\tau\right]=\frac{c(t)}{\kappa(t)},  \qquad h(t)=\frac{3}{4}.
\end{equation}
Imposing $h=3/4$, we obtain the integral equation
\begin{equation}\label{relat-2}
 4q \kappa(t) \exp\left[\int_0^t (1-2c(\tau))d\tau\right]=c(t).
\end{equation}

\begin{thm}
Let us suppose that $h=\frac{3}{4}$ and $c(t),\kappa(t)$ are
related by \eqref{relat-2} with a given $q>0$. Then, system
\eqref{model_eq2} is solvable.
\end{thm}

The simplest way to find explicit examples is to prescribe $c(t)$
and use \eqref{relat-2} to find the adequate $\kappa(t)$. For instance, if $c$ is constant, then
\eqref{relat-2} gives
$$
\kappa(t)=\frac{c}{4q}\exp\left[(2c-1)t\right].
$$
If $c=1/2$, then $\kappa$ is an arbitrary positive constant. This
case $c=1/2,h=3/4$ was identified in \cite[Subcase
(B.1)]{NucciSanchini2015}. We are generalizing this particular
case to an exponentially increasing (if $c>1/2$) or decreasing (if
$c<1/2$) carrying capacity.

It is also possible to identify periodic coefficients. For
instance, if $c(t)=\frac{1}{2}+\lambda \cos t$, then
$$
\kappa(t)=\left(\frac{1}{2}+\lambda \cos
t\right)\frac{1}{4q}\exp\left[2\lambda \sin t\right].
$$

Conversely, it is possible to consider a prescribed carrying
capacity $\kappa(t)$ and look for a suitable $c(t)$. To do it, one observes that a derivation of \eqref{relat-2} gives
$$
\dot c(t)=\left(1-\frac{\dot\kappa(t)}{\kappa(t)}\right)c-2c^2,
$$
which can be identified as a Bernouilli's equation and hence explicitly solvable. The correct initial condition comes from evaluating \eqref{relat-2} in $t=0$, hence $c(0)=4q\kappa(0)$.

As an example, we perform a detailed analysis or the case of a constant carrying capacity $\kappa$. First, evaluating \eqref{relat-2} in $t=0$, one realizes
that the correct initial condition is $c(0)=4q \kappa$. Then, a
simple derivation on $t$ leads to the separable equation
$$
\dot c=c(1-2c),
$$
together with the initial condition $c(0)=4q \kappa$ has the
explicit solution
\begin{equation}\label{c-form}
c(t)=\frac{1}{C_0 e^{-t}+2}
\end{equation}
with $C_0=\displaystyle\frac{1-8q \kappa}{4q \kappa}$.

Hence, the new form of our initial equation boils down to
\begin{equation}\label{2nd-model-integ}
  \ddot{y}+r(t)\dot{y}+\alpha(t)y=\frac{5}{4}\frac{\dot{y}^2}{y}+\exp\left[-2\int_0^t r(\tau)d\tau\right]y^2.
\end{equation}
with coefficients
\begin{equation}\label{coeff}
  r(t)=c(t)-\frac{1}{2}, \quad \alpha(t)=\frac{3}{4}-c(t).
\end{equation}
Provided $c(t)$ is of the form \eqref{c-form}, the general
solution is given by
\begin{equation}\label{gen-sol}
  y(t)=\exp\left[2\int_0^t r(\tau)d\tau\right](Au^2+2Buv+Cv^2)^{-2},  \quad (AC-B^2)W^2=q,
\end{equation}
where $u, v$ are two linearly independent solutions of the base
equation \eqref{base} with $p$ replaced by $-\alpha(t)/4$.

A simple computation gives
\begin{equation}\label{r}
\exp\left[2\int_0^t r(\tau)d\tau\right]=\exp\left[\ln (C_0+2 e^t)-t-\ln
(C_0+2)\right]=\frac{C_0e^{-t}+2}{C_0+2}.
\end{equation}
On the other hand, the base equation \eqref{base} is written
explicitly as
\begin{equation}\label{base-2}
\ddot z+I(t)z=0
\end{equation} with
$$
I(t)=-\frac{3}{16}+\frac{1}{4(C_0 e^{-t}+2)}-\frac{C_0
e^{-t}(C_0e^{-t}+8)}{16(C_0 e^{-t}+2)^2}.
$$
Observe that $I(t)\to -1/16$ when $t\to +\infty$, therefore
\eqref{base-2} has an exponential dichotomy and $u,v$ can be taken
such that
$$
\lim_{t\to+\infty}u(t)=+\infty, \quad  \lim_{t\to +\infty}v(t)=0.
$$
In view of \eqref{r} and \eqref{gen-sol}, this implies that
generically (i.e., for most of the solutions), $y(t)\to 0$ when
$t\to +\infty$, that is, asymptotic extinction for $y$.

For the particular case of $C_0=0$, then \eqref{base-2} is just
$\ddot z-\frac{1}{16}z=0$ and we can take $u(t)=e^{t/4},
v(t)=e^{-t/4}$, for which the Wronskian is $W=-1/2$. Then,
\eqref{gen-sol} gives
$$
y(t)=2(Ae^{t/2}+2B+Ce^{-t/2})^{-2},\quad AC-B^2=4q,
$$
Now, $x(t)$ is obtained from \eqref{x}.

More surprisingly, when $C_0\ne 0$ equation \eqref{base-2} is
still integrable, with fundamental system
$$
u=e^{-t/2}(C_0+2 e^t)^{3/4}, \qquad v=e^{-t/2}(C_0+2 e^t)^{1/4},  \quad W=-1.
$$
The solution $y(t)$ can be expressed as
$$y(t)=\frac{e^t}{C_0+2}[A(C_0+2 e^t)+2B(C_0+2 e^t)^{1/2}+C]^{-2},  \quad  AC-B^2=q.$$

\section*{Acknowledgments}
F.G. would like to thank Prof. P.J. Torres for the kind invitation and support to visit the Department of Applied Mathematics of the University of Granada. This work is partially supported by MICENO and ERDF project MTM2017-82348-C2-1-P.

\bibliographystyle{unsrt}

\end{document}